\documentclass[12pt, psamsfonts]{amsart}

\usepackage{thmtools}

\usepackage{centernot}

\usepackage{tikz}
\usetikzlibrary{arrows.meta}

\usepackage[shortlabels]{enumitem}

\usepackage[spacing=true,kerning=true,babel=true,tracking=true]{microtype}

\usepackage[alphabetic, nobysame]{amsrefs}

\usepackage{float}

\usepackage{hyperref}
\hypersetup{
	pdfstartview={XYZ null null 1.00}, 
	pdfpagemode=UseNone, 
	colorlinks,
	breaklinks, 
	linkcolor=blue,
	urlcolor=blue, 
	anchorcolor=blue,
	citecolor=blue
}

\usepackage[capitalise]{cleveref}
\crefformat{footnote}{#2\footnotemark[#1]#3}
\usepackage{nameref}

\usepackage{geometry}
\usepackage{comment}

\usepackage{xspace}
\usepackage{ bbold }
\usepackage{graphicx}

\usepackage{etoolbox}

\usepackage{thmtools}
\usepackage{thm-restate}
\declaretheorem[name=Theorem,numberwithin=section]{thm}

\usepackage{amssymb}

\makeatletter
\def\subsection{\@startsection{subsection}{2}%
  \z@{.8\linespacing\@plus.7\linespacing}{.5\linespacing}%
  {\normalfont\bfseries}}

\def\@seccntformat#1{%
  \protect\textup{%
    \protect\@secnumfont
    \expandafter\protect\csname format#1\endcsname % <--- added
    \csname the#1\endcsname
    \protect\@secnumpunct
  }%
}
\patchcmd{\@startsection}
  {\@afterindenttrue}
  {\@afterindentfalse}
  {}{}
 \makeatother

\makeatletter
\def\l@section{\@tocline{1}{5pt}{0pc}{}{}}
\renewcommand{\tocsection}[3]{%
	\indentlabel{\@ifnotempty{#2}{\makebox[20pt][l]{%
				\ignorespaces#1 #2.\hfill}}}\sc #3\dotfill}

\newdimen{\tocsubsecmarg}
\def\l@subsection{\@tocline{2}{3pt}{0pc}{\tocsubsecmarg}{}}
\renewcommand{\tocsubsection}[3]{%
	\indentlabel{\@ifnotempty{#2}{\makebox[30pt][l]{%
				\ignorespaces#1 #2.\hfill}}}#3\dotfill}
\makeatother
\let\oldtocsubsection=\tocsubsection
\renewcommand{\tocsubsection}[2]{\hspace{3em} \oldtocsubsection{#1}{#2}}

\numberwithin{equation}{section}

\theoremstyle{plain}

\crefname{prop}{Proposition}{Propositions}

\crefname{obs}{Observation}{Observations}

\crefname{cor}{Corollary}{Corollaries}
\newtheorem{cor}[equation]{Corollary}

\makeatletter
\def\@empty{}
\def\ifemptycredit#1{%
	\def\tmp{#1}%
	\ifx\tmp\@empty%
	\else%
	{~(#1)}%
	\fi%
}
\makeatother

\newenvironment{namedthm*}[2][]{
\medskip\par\noindent \textbf{#2}\ifemptycredit{#1}\textbf{.}\itshape\xspace
}{\medskip}

\theoremstyle{definition}

\crefname{defn}{Definition}{Definitions}

\crefname{example}{Example}{Examples}

\crefname{question}{Question}{Question}

\crefname{ergprop}{Nonsingular Pointwise Ergodic Property}{Ergodic Properties}

\theoremstyle{remark}

\newtheorem*{acknowledge}{Acknowledgments}

\crefname{remark}{Remark}{Remarks}

\crefname{claim}{Claim}{Claims}
\newtheorem{claim}[equation]{Claim}
\newtheorem*{claim*}{Claim}

\crefname{assumption}{Assumption}{Assumptions}

\declaretheoremstyle[
spaceabove=\topsep, 
spacebelow=6pt,
headfont=\normalfont\itshape,
notefont=\normalfont, notebraces={(}{)},
bodyfont=\normalfont,
postheadspace=4pt,
qed=\mbox{\smaller[4]$\boxtimes$}
]{claimproofstyle}
\declaretheorem[name={Proof of Claim}, style=claimproofstyle, unnumbered]{pf}

\let\#\mathbb
\let\@\mathcal

\newcommand{\N}{\mathbb{N}}

\newcommand{\e}{\varepsilon}

\renewcommand{\phi}{\varphi}

\newcommand{\actson}{\curvearrowright}

\newcommand{\ie}{i.e.,\xspace}

\newcommand*{\defeq}{\mathrel{\vcenter{\baselineskip0.5ex \lineskiplimit0pt \hbox{\scriptsize.}\hbox{\scriptsize.}}}=}

\newcommand{\set}[1]{\left\{ #1 \right\}}

\usepackage{color}

\definecolor{purple}{RGB}{116,0,159}

\definecolor{vert}{RGB}{7,126,26}

\newcommand{\Act}{\text{Act}}
\newcommand{\Aut}{\text{Aut}}

\newcommand{\symdif}{\mathbin{\triangle}}

\let\@\mathcal
\let\~\widetilde

\title{Every countable group has the \\ nonsingular weak Rokhlin property}

\author{Garrett Ervin}
\address{(Garrett Ervin, Bogl\'arka Geh\'er, Anett Imolay-Kocsis, and M\'at\'e Pálfy) Institute of Mathematics, E\"otv\"os Lor\'and University, Budapest, Hungary}
\email{garrette@ttk.elte.hu}
\author{Boglárka Gehér}
\email{bogigeher@gmail.com}
\author{Anett Imolay-Kocsis}
\email{sakkboszi@gmail.com}
\author{Máté Pálfy}
\email{palfymateandras@gmail.com}
\author{Jenna Zomback}
\address{(Jenna Zomback) Mathematics, Amherst College, Amherst, MA, USA}
\email{jzomback@amherst.edu}

\date{\today}

\begin{document}

\begin{abstract}
    We prove that for every countable group $\Gamma$, in the space of nonsingular actions by $\Gamma$, there is an action whose orbit under conjugation is dense. This is the nonsingular analogue of the result of Glasner, Thouvenot, and Weiss that every countable group $\Gamma$ is weak Rokhlin, \ie there is a dense orbit in the space of probability measure preserving $\Gamma$-actions. 
\end{abstract}

\maketitle

\section{Introduction}

A countably infinite group $\Gamma$ has the \textit{weak Rokhlin property} if there is a probability measure preserving (pmp) action $\alpha$ by $\Gamma$ on a standard probability space $(X, \mu)$ whose conjugates by pmp automorphisms of $X$ are dense in the Polish space of pmp actions $\Act(\Gamma, X, \mu)$. 

The following result was proved by Glasner, Thouvenot, and Weiss in \cite{GTW}, and independently by Hjorth in unpublished work; see \cite{Kechris}*{Theorem 10.7}. 

\begin{thm}[Glasner-Thouvenot-Weiss, Hjorth]\label{weak_rokhlin_pmp}
    Every countably infinite group has the weak Rokhlin property.
\end{thm}

The desired action $\alpha$ is obtained by taking the diagonal action on $X^{\N}$ by a sequence of actions $(\alpha_i)_{i \in \N}$ that is dense in $\Act(\Gamma, X, \mu)$, and then conjugating $\alpha$ by an isomorphism $\varphi: X \rightarrow X^{\N}$. 

We say that a countably infinite group $\Gamma$ has the \textit{nonsingular weak Rokhlin property} if there is a measure-class preserving (mcp) action $\Gamma \actson X$ whose conjugates by mcp automorphisms of $X$ are dense in the space of mcp actions $\Act^*(\Gamma, X, \mu)$. 

Since a countable product of mcp actions is not necessarily mcp, mimicking the proof from \cite{GTW} does not yield an mcp version of Theorem \ref{weak_rokhlin_pmp}. 
Nevertheless, we show that the analogue of Theorem \ref{weak_rokhlin_pmp} holds in the mcp setting.

We build the desired action by concatenating copies of the actions from a dense sequence $(\alpha_i)_{i \in \N}$ over a countable interval partition of $(X, \mu) = ([0, 1], \lambda)$.  

\begin{thm}\label{weak_rokhlin}
Every countably infinite group has the nonsingular weak Rokhlin property.
\end{thm}

\begin{acknowledge}
This work was initiated at the \textit{Group actions and equivalence relations} workshop in Budapest, Hungary in July 2026, which was supported by HAS Momentum Grant 2022-58, NKFIH KKP 139502.
The authors thank the organizers Miklós Abért, Clark Lyons, Petr Naryshkin, Grigory Terlov, and Zoltán Vidnyánszky for inviting us and for creating a warm and collaborative environment.

The first author is supported by the Momentum MSCA Programme, which is co-funded by the European Commission through the HORIZON-MSCA-2023-COFUND programme and the Secretariat of the Hungarian Academy of Sciences (MTA).

The fourth author is supported by the National Research, Development and Innovation Office -- NKFIH, grants no.~152822. 
\end{acknowledge}

\noindent \textit{AI use statement}. No AI tools were used in the research or writing of this paper. 

\section{Preliminaries}

Let $(X, \mu)$ be a standard atomless probability space. A measurable bijection $T: X \rightarrow X$ is called a \textit{measure-class preserving} (or \textit{mcp}, or \textit{nonsingular}) \textit{transformation} if for every null set $A \subseteq X$, $T(A)$ is null. 

Let $\Aut^*(X, \mu)$ denote the group of mcp transformations of $X$ equipped with the weak topology, where two transformations are identified when they are equal modulo a null set. Under this topology, a sequence of transformations $T_n$ converges to $T$ if
\[
    \mu(T_n(A) \symdif T(A))\to 0
\]
for every measurable set $A \subseteq X$, and
\[
\frac{d(\mu \circ T_n)}{d\mu} \to \frac{d(\mu \circ T)}{d\mu}
\]
in the $L_1$ norm. This topology makes $\Aut^*(X, \mu)$ a Polish group; see \cite{DS}*{Section 4} for a reference.

Fix a countably infinite discrete group $\Gamma$. Let $\Act^*(\Gamma, X, \mu)$ denote the space of group actions by $\Gamma$ on $X$ by mcp transformations, viewed as a subspace of the product space $(\Aut^*(X, \mu))^{\Gamma}$: 
\[
\Act^*(\Gamma, X,\mu) = \set{\alpha \in (\Aut^*(X,\mu))^\Gamma \; : \; \alpha \text{ is a group action}}.
\]

Given an action $\alpha \in \Act^*(\Gamma, X, \mu)$ and group element $\gamma \in \Gamma$, we write $\alpha^{\gamma}$ for the mcp transformation $\alpha(\gamma) \in \Aut^*(X, \mu)$. 

Observe that $\Act^*(\Gamma, X, \mu)$ is a closed subspace of the Polish space $(\Aut^*(X, \mu))^{\Gamma}$ and hence is itself Polish. 
The group $\Aut^*(X, \mu)$ acts naturally on the space $\Act^*(\Gamma, X, \mu)$ by conjugation. 

\section{Proof of \cref{weak_rokhlin}}

Theorem \ref{weak_rokhlin} asserts that there is always an action $\alpha \in \Act^*(\Gamma, X, \mu)$ whose orbit under conjugation by $\Aut^*(X, \mu)$ is dense in $\Act^*(\Gamma, X, \mu)$.
We first sketch the idea of the proof. 

\medskip 

\noindent
\textit{Proof sketch:}
Let $(\alpha_n)_{n \in \N}$ be dense in $\Act^*(\Gamma,[0,1],\lambda)$. 
We construct $\Gamma\actson^\alpha ([0,1],\lambda)$ by partitioning $[0,1]$ into $\N$-many nondegenerate closed intervals (modulo endpoints) $(I_n)_{n \in \N}$, and we let $\alpha$ act on $I_n$ by $\alpha_n$.
We then fix $n \in \N$ and approximate $\alpha_n$ by $S \alpha S^{-1}$, where $S$ maps $I_n$ to the interval $[0,1-\e]$.
The idea is that since $S \alpha S^{-1}$ and $\alpha_n$ more or less agree on $1-\e$ fraction of the space, this is a good approximation.

 \begin{proof}
Let $\Gamma$ be a countable group.
Since $([0,1], \lambda)$ and $(X , \mu)$ are Borel-isomorphic measure spaces, we may assume $(X, \mu) = ([0,1], \lambda)$.
 \subsubsection*{Construction of $\alpha$}
    Let $(\alpha_n)_{n \in \N}$ be a dense sequence in $\Act^*(\Gamma,[0,1],\lambda)$. 
    It suffices to construct an action $\alpha$ whose conjugacy class contains each $\alpha_n$ in its closure.
    Let
    \[
        [0,1]=\bigsqcup_{n\in \N} I_n
    \]
    be a partition, modulo endpoints, of $[0, 1]$ into nondegenerate closed intervals. For each $n \in \N$, label the endpoints of $I_n$ by writing $I_n = [a_n, b_n]$, and let
    \[
        \iota_n\colon [0,1]\longrightarrow I_n
    \]
    be the affine bijection defined by 
    \[
    \iota_n(x) = (b_n - a_n)x + a_n.
    \]
    Observe that each $\iota_n$ is measure-class preserving.  

    Define an action
    $\alpha\in\Act^*(\Gamma,[0,1],\lambda)$ by
    \[
        \alpha^\gamma(x)
        =
        \iota_n\circ\alpha_n^\gamma\circ\iota_n^{-1}(x)
        \qquad
        \text{for }x\in I_n.
    \]
    Thus each interval $I_n$ is $\alpha$-invariant, and the restriction of $\alpha$ to $I_n$ is isomorphic to $\alpha_n$. 
    
    Since each $\alpha_n^\gamma$ is mcp, so is each $\alpha^\gamma$. Indeed, if $A \subseteq [0,1]$ is measurable, then $A$ is null if and only if $A \cap I_n$ is null for every $n \in \N$. Since $\alpha_n^\gamma$ is mcp, $A \cap I_n$ is null exactly when $ \iota_n\circ\alpha_n^\gamma\circ\iota_n^{-1}(A \cap I_n) = \alpha^\gamma(A \cap I_n)$ is null.
    
    %A null if A cap I_n null for all n (ctbl additivity) iff ... is null for all n (using that X is null if affine scaling of X is null) iff alpha^gamma(A cap I_n) is null for all n iff alpha^gamma(A) is null (using that measure of alpha^gamma(A) = sum of measures of alpha^gamma(A cap I_n)

\subsubsection*{Approximation of $\alpha_k$}
    Fix $k\in \N$. We show that $\alpha_k$ belongs to the closure of the
    conjugacy class of $\alpha$. For $0<\varepsilon<1$, put
    \[
        J_\varepsilon=[0,1-\varepsilon],
        \qquad
        K_\varepsilon=(1-\varepsilon,1],
    \]
    and define
    \[
        h_\varepsilon\colon [0,1]\longrightarrow J_\varepsilon,
        \qquad
        h_\varepsilon(x)=(1-\varepsilon)x.
    \]
    Choose a measure-class preserving Borel isomorphism
    \[
        \tau_\varepsilon\colon
        K_\varepsilon\longrightarrow [0,1]\setminus I_k,
    \]
    and define $\pi_\varepsilon\in\Aut^*([0,1],\lambda)$ by
    \[
        \pi_\varepsilon(x)
        =
        \begin{cases}
            \iota_k\circ h_\varepsilon^{-1}(x),
                & x\in J_\varepsilon,\\[2mm]
            \tau_\varepsilon(x),
                & x\in K_\varepsilon.
        \end{cases}
    \]
    In words, $\pi_\varepsilon$ maps the set $J_\e$ onto the interval $I_k$, and maps $K_\e$ onto the union of all other intervals.

    Consider the conjugate action
    \[
        \beta_\varepsilon
        =
        \pi_\varepsilon^{-1}\circ
        \alpha \circ
        \pi_\varepsilon.
    \]

        For $x\in J_\varepsilon$, we have
    \begin{equation}\label{e.conjugation}
        \beta_\varepsilon^\gamma(x)
        =
        \pi_\varepsilon^{-1}
        \circ\iota_k
        \circ\alpha_k^\gamma
        \circ\iota_k^{-1}
        \circ\pi_\varepsilon(x)
        =
        h_\varepsilon
        \circ\alpha_k^\gamma
        \circ h_\varepsilon^{-1}(x).
        \tag{$\star$}
    \end{equation}
    Moreover, $K_\varepsilon$ is $\beta_\varepsilon$-invariant.
        Intuitively, $\beta_\varepsilon$ approximates $\alpha_k$ by acting by a slightly compressed version of $\alpha_k$ on $J_\e$ and compressing the rest of $\alpha$ to $K_\e$. 
        
    \subsubsection*{Proof of convergence}
    We claim that in the convergence topology on $\Act^*(\Gamma, [0,1],\lambda)$,
    \[
        \beta_\varepsilon\longrightarrow\alpha_k
        \qquad\text{as }\varepsilon\longrightarrow0.
    \]
    
    By definition of the topology on $\Act^*(\Gamma,[0,1],\lambda)$, it suffices to show that for a fixed measurable $B \subseteq [0, 1]$ and $\gamma \in \Gamma$, we have
  
    \begin{enumerate}
        \item $\lim_{\e \to 0} \lambda(\beta_\varepsilon^\gamma(B) \symdif \alpha_k^\gamma(B)) = 0$, and
        \item $\lim_{\e \to 0} ||\frac{d (\alpha_k^\gamma \circ \lambda)}{d \lambda} - \frac{d (\beta_\varepsilon^\gamma \circ \lambda)}{d\lambda}||_1 = 0$.
    \end{enumerate}

    We establish (1) and (2) as Claims \ref{firstclaim} and \ref{secondclaim} below. 

    Fix $B$ and $\gamma$. For convenience, we abbreviate
    \[
        T \defeq \alpha_k^\gamma
        \; \text{ and } \;
        S_\varepsilon \defeq \beta_\varepsilon^\gamma.
    \]

    \begin{claim} \label{firstclaim}
    $\lim_{\e \to 0} \lambda(S_\e(B) \symdif T(B)) = 0$.
    \end{claim}

    \begin{pf}
    For a measurable subset $A \subseteq [0,1]$, define
    \[
        A^{\varepsilon}
        \defeq
        h_\varepsilon^{-1}(A\cap J_\varepsilon).
    \]

    Observe that if $C$ is a finite union of intervals, then
    \[
        \lambda(C^{\varepsilon}\mathbin{\triangle}C)\longrightarrow0
        \text{ and }
        \lambda\bigl(h_\varepsilon(C)\mathbin{\triangle} C \bigr)
        \longrightarrow0.
    \]
    
    Since $B$ is measurable, there is a finite union of intervals $C$ making $\lambda(B \mathbin{\triangle} C)$ arbitrarily small. Then $\lambda(B^{\varepsilon} \mathbin{\triangle} C^{\varepsilon}) = \lambda((B \mathbin{\triangle} C)^{\varepsilon})$ and $\lambda(h_{\varepsilon}(B) \mathbin{\triangle} h_{\varepsilon}(C)) = \lambda(h_{\varepsilon}(B \mathbin{\triangle} C))$ can be made arbitrarily small by letting $\varepsilon$ approach $0$. Since
    \[
    (B^{\varepsilon} \mathbin{\triangle} B) \subseteq (B^{\varepsilon} \mathbin{\triangle} C^{\varepsilon}) \cup (C^{\varepsilon} \mathbin{\triangle} C) \cup (C \mathbin{\triangle} B),
    \]
    it now follows that 
    \[
        \lambda(B^{\varepsilon}\mathbin{\triangle}B)\longrightarrow0,
    \]
    and a similar calculation shows
    \[
        \lambda\bigl(h_\varepsilon(B)\mathbin{\triangle} B \bigr)
        \longrightarrow0.
    \]

    Since $T$ is measure-class preserving, the measures $\lambda$ and $\lambda \circ T$ are equivalent. 
    Consequently,
    \[
        \lambda\bigl(T(B^{\varepsilon})
        \mathbin{\triangle}T(B)\bigr)
        =
        \lambda\bigl(T(B^{\varepsilon}\mathbin{\triangle}B)\bigr)
        \longrightarrow0.
    \]
    Observe that by \ref{e.conjugation}, we have
    $S_\varepsilon(B\cap J_\varepsilon)
        =
        h_\varepsilon\bigl(T(B^{\varepsilon})\bigr)$
    and $S_\varepsilon(B\cap K_\varepsilon)\subseteq K_\varepsilon$.
    Therefore,
    \begin{align*}
        \lambda\bigl(S_\varepsilon(B)
        \mathbin{\triangle}T(B)\bigr)
        &= \lambda\bigl(S_\varepsilon((B \cap J_\e)\sqcup S_\varepsilon(B \cap K_\e))
        \mathbin{\triangle}T(B)\bigr)
        \\
        &= \lambda\bigl((h_\e(T(B^\e)) \sqcup S_\e(B \cap K_\e))
        \mathbin{\triangle}T(B)\bigr)
        \\
        &\leq \lambda\bigl(h_\e(T(B^\e))
        \mathbin{\triangle}T(B)\bigr)
        +
        \lambda(S_\e(B \cap K_\e))
        \\
        &\leq
        \lambda\bigl(
            h_\varepsilon(T(B^{\varepsilon}))
            \mathbin{\triangle}
            h_\varepsilon(T(B))
        \bigr)
        +
        \lambda\bigl(
            h_\varepsilon(T(B))
            \mathbin{\triangle}
            T(B)
        \bigr)
        +\varepsilon\\
        &=
        (1-\varepsilon)
        \lambda\bigl(
            T(B^{\varepsilon})
            \mathbin{\triangle}
            T(B)
        \bigr)
        +
        \lambda\bigl(
            h_\varepsilon(T(B))
            \mathbin{\triangle}
            T(B)
        \bigr)
        +\varepsilon,
    \end{align*}
    and the right-hand side tends to zero.
    \qedhere
    \end{pf}

    It remains to verify (2). 
    We abbreviate
    \[
        \rho_T
        =
        \frac{d (T \circ \lambda)}{d \lambda}
        \; \text{ and } \;
        \rho_{S_\varepsilon}
        =
        \frac{d (S_\varepsilon \circ \lambda)}{d \lambda}.
    \]

    \begin{claim} \label{secondclaim}
        $\lim_{\e \to 0}\|\rho_{S_\varepsilon}-\rho_T\|_1 = 0$.
    \end{claim}

    \begin{pf}
    Since $h_\varepsilon$ rescales Lebesgue measure by a constant factor, this
    factor cancels under conjugation, and hence
    \[
        \rho_{S_\varepsilon}(x)
        =
        \rho_T\bigl(h_\varepsilon^{-1}(x)\bigr)
        \qquad
        \text{for a.e. }x\in J_\varepsilon.
    \]
    Furthermore, since $K_\varepsilon$ is $S_\varepsilon$-invariant,
    \[
        \int_{K_\varepsilon}
        \rho_{S_\varepsilon}d\lambda
        =
        \lambda(K_\varepsilon)
        =
        \varepsilon.
    \]
    It follows that
    \begin{align*}
        \|\rho_{S_\varepsilon}-\rho_T\|_1
        &\leq
        \int_{J_\varepsilon}
        \left|
            \rho_T\bigl(h_\varepsilon^{-1}(x)\bigr)-\rho_T(x)
        \right| \; 
        d\lambda(x)
        +
        \int_{K_\varepsilon}\rho_{S_\varepsilon} \; d\lambda
        +
        \int_{K_\varepsilon}\rho_T \; d\lambda\\
        &=
        \int_{J_\varepsilon}
        \left|
            \rho_T\bigl(h_\varepsilon^{-1}(x)\bigr)-\rho_T(x)
        \right|
         \; d \lambda(x)
        +
        \varepsilon
        +
        \int_{K_\varepsilon}\rho_T \; d \lambda.
    \end{align*}
    After approximating with a simple function, it is easy to see that the first term tends to zero, while the last term tends to zero by the absolute
    continuity of the integral. Thus
    \[
        \|\rho_{S_\varepsilon}-\rho_T\|_1\longrightarrow0.
        \qedhere
    \] 
    \end{pf}

     We now have that $\beta_\varepsilon\longrightarrow\alpha_k$, as desired.
\end{proof}

Observe also that $\alpha$ is free exactly when each of the actions $\alpha_k$ is free, and as long as \emph{one} of the $\alpha_k$ is free, $\alpha$ is \emph{somewhere free}, \ie free on an invariant positive measure set.
We also note that regardless of the sequence $(\alpha_i)_{i \in \N}$, $\alpha$ is not ergodic since each subinterval is invariant.
We thus obtain the following.

\begin{cor}
    For any countable group $\Gamma$ and standard probability space $(X,\mu)$, the following sets are dense in $\Act^*(\Gamma, X, \mu)$.
    \begin{enumerate}[(a)]
        \item the set of non-ergodic actions
        \item the set of non-free actions
        \item the set of somewhere-free actions
    \end{enumerate}
\end{cor}

\bibliographystyle{amsalpha} 
\bibliography{mcp} 
\end{document}